\documentclass[12pt]{article}
\usepackage{amssymb}
\usepackage{latexsym,bm}
\usepackage{graphicx}
\usepackage{amsmath}
\usepackage{amsthm}
\usepackage{mathrsfs}

\newtheorem{thm}{Theorem}
\newtheorem{proof of theorem }{Proof of Theorem }

\newtheorem{nota}{Notation}
\newtheorem{defi}{Definition}

\newtheorem{lem}[thm]{Lemma}
\newtheorem{cor}[thm]{Corollary}

\newcommand{\bea}{\begin{eqnarray*}}
\newcommand{\eea}{\end{eqnarray*}}
\newcommand{\be}{\begin{equation}}
\newcommand{\ee}{\end{equation}}
\newcommand{\ben}{\begin{eqnarray*}}
\newcommand{\een}{\end{eqnarray*}}

\date{}
\begin{document}
\title{The maximum size of a nonhamiltonian-connected graph with  given order and minimum degree\footnote{E-mail addresses:
{\tt mathdzhang@163.com}.}}
\author{\hskip -10mm Leilei Zhang\\
{\hskip -10mm \small School of Mathematics and Statistics, Central China Normal University, Wuhan 430079, China}}\maketitle
\begin{abstract}
In this paper, we  determine the maximum size of a nonhamiltonian-connected graph with prescribed order and minimum degree. We also characterize the extremal graphs that attain this maximum size. This work generalizes a previous result obtained by Ore [ J. Math. Pures Appl. 42 (1963) 21-27] and further extends a theorem proved by Ho, Lin, Tan, Hsu, and Hsu [Appl. Math. Lett. 23 (2010) 26-29]. As a corollary of our main result, we determine the maximum size of a $k$-connected nonhamiltonian-connected graph with a given order.
\end{abstract}
{\bf Key words.} Nonhamiltonian-connected; minimum degree; maximum size.

\noindent{\bf Mathematics Subject Classification.} 05C30, 05C35, 05C45

\section{ Introduction}
The graphs under consideration are simple and undirected. We use the book of West \cite{9} for common terminology and notation not defined here. The {\it order} of a graph is its number of vertices, and the {\it size} its number of edges. Denote by $V(G)$ and $E(G)$ the vertex set and edge set of a graph $G$ respectively. Let $e(G)$ denote the size of $G.$ For two vertices $u$ and $v$, the symbol $u\leftrightarrow v$ means that $u$ and $v$ are {\it adjacent} and $u\nleftrightarrow v$ means that $u$ and $v$ are {\it nonadjacent}. Let $d_G(v)$ and $N_G(v)$ denote, respectively, the {\it degree} and {\it neighbourhood} of a vertex $v\in V (G)$.  We denote by $\delta(G)$ the {\it minimum degree} of a graph $G.$ For a subgraph $H$ of $G$, $N_H(v)\triangleq N_G(v)\cap V(H)$ and the degree of $v$ in $H$ is $d_H(v)\triangleq |N_G(v)\cap V(H)|$ .

For graphs we will use equality up to {\it isomorphism}, so $G_1=G_2$ means that $G_1$ and $G_2$ are isomorphic. $\overline{G}$ denotes the complement of a graph $G.$ For two graphs $G$ and $H,$ $G\vee H$ denotes the {\it join} of $G$ and $H,$ which is obtained from the {\it disjoint union} $G+H$ by adding edges joining every vertex of $G$ to every vertex of $H.$ $K_n$ denotes the complete graph of order $n.$ A cycle (path) containing all vertices of a graph $G$ is called a {\it Hamilton cycle (path)} of $G$. A graph $G$ is {\it hamiltonian} if it contains a Hamilton cycle; otherwise it is a nonhamiltonian graph. A {\it pancyclic} graph is a graph that contains cycles of all possible lengths from three up to the number of vertices in the graph. A graph $G$ is {\it hamiltonian-connected} if there exists a hamiltonian path joining any two vertices of $G$.

The problem of determining whether a given graph is Hamiltonian (or traceable, or hamiltonian-connected) is NP-complete \cite{11}. There is a vast literature \cite{12,19,14,15,13} in graph theory devoted to obtaining sufficient conditions for Hamiltonicity. In 1961, Ore \cite{8} determined the maximum size of a nonhamiltonian graph with a given order and also determined the extremal graphs.

\begin{thm}{\rm (Ore \cite{8})}
The maximum size of a nonhamiltonian graph of order $n\ge 3$ is $\binom{n-1}{2}+1$ and this size is attained by a graph $G$ if and only if $G=K_1\vee (K_{n-2}+K_1)$ or $G=K_2\vee \overline{K_3}.$
\end{thm}

Bondy \cite{10} gave a new proof of Theorem 1. Obviously, if the size of a graph exceeds $\binom{n-1}{2}+1$, then the graph is hamiltonian. It is natural to ask the same question by putting constraints on the graphs. In 1962 Erd\H{o}s \cite{6} determined the maximum size of a nonhamiltonian graph of order $n$ and minimum degree at least $k.$

\begin{thm}{\rm (Erd\H{o}s \cite{6})}
Let $n$, $\delta$ be integers with $1\le \delta\le \lfloor \frac{n-1}{2}\rfloor.$  If $G$ is a nonhamiltonian graph of order $n$ with minimum degree at least $\delta$, then
$$
e(G)\le {\rm max}\left\{\binom{n-\delta}{2}+ \delta^2,\, \binom{n-\left\lfloor\frac{n-1}{2}\right\rfloor}{2}+ \left\lfloor\frac{n-1}{2}\right\rfloor^2\right\}.
$$
\end{thm}

An $s$-clique is a clique of cardinality $s$. Motivated by Erd\H{o}s' work, F\"{u}redi, Kostochka and Luo \cite{3} determined the maximum number of $s$-cliques in nonhamiltonian graphs. Recently, Zhang \cite{16} determined the maximum size of a nonhamiltonian graph with given order and minimum degree.

\begin{thm}{\rm (Zhang \cite{16})}
Given positive integers $n$ and $\delta$ with $n\ge \delta+1,$ let $t=\left\lfloor\frac{n-1}{2}\right\rfloor.$ Then the maximum size of a nonhamiltonian $2$-connected graph of order $n$ and minimum degree $\delta$ is
$$
{\rm max}\left\{\binom{n-\delta}{2}+ \delta^2, \binom{n-t}{2}+t(t-1)+\delta\right\}.
$$
\end{thm}

The graphs that attain $\phi(n,\delta)$ are also characterized in \cite{16}. Bondy suggested that almost all non-trivial conditions for Hamiltonicity of a graph should also imply pancyclicity. Inspired by this, Schmeichel and Hakimi \cite{1} proved the following result.

\begin{thm}{\rm(Schmeichel and  Hakimi \cite{1})}
Let $\psi(n,\delta)$ denote the maximum size of a nonpancyclic graph of order $n$ and minimum degree $\delta.$ Then
$$
\psi(n,\delta)=\begin{cases}
\binom{n-\delta}{2}+\delta^2\quad {\rm if}\,\,n\,\,\,{\rm is}\,\,{\rm odd}\,\,{\rm and}\,\,1\le \delta \le \frac{n+5}{6}\,\,{\rm or}\,\,
n\,\,{\rm is}\,\,{\rm even}\,\,{\rm and}\,\,1\le \delta \le \frac{n+8}{6},\\
\frac{3n^2-8n+5}{8}+\delta\quad {\rm if}\,\,n\,\,{\rm is}\,\,{\rm odd}\,\,{\rm and}\,\,\frac{n+5}{6}\le \delta\le\frac{n-1}{2},\\
\frac{3n^2-10n+16}{8}+\delta\quad {\rm if}\,\,n\,\,{\rm is}\,\,{\rm even}\,\,{\rm and}\,\,\frac{n+8}{6}\le \delta\le\frac{n-2}{2}.
\end{cases}
$$
\end{thm}

Hamiltonian-connected graph is an important type of hamiltonian graphs. There has been extensive research on the properties of hamiltonian-connected graphs, see \cite{20,18}. Some sufficient conditions for a graph to be hamiltonian-connected have also been found in \cite{4,2}. One cornerstone in this direction is the following celebrated Ore's theorem.

\begin{thm}{\rm (Ore \cite{2})}
The maximum size of a nonhamiltonian-connected graph of order $n$ is $\binom{n-1}{2}+2$ and this size is attained by a graph $G$
if and only if $G=K_3\vee \overline{K_3}$ or $G= K_2\vee( K_{n-3}+ K_1).$
\end{thm}

Note that if the number of edges of a graph exceeds $\binom{n-1}{2}+2$, then the graph is hamiltonian-connected. By imposing minimum degree as a new parameter, Ho, Lin, Tan, Hsu and Hsu \cite{4} determined the maximum size of a nonhamiltonian-connected graph with given order and minimum degree at least $\delta.$ Obviously, any hamiltonian-connected graph of order at least $4$ is $3$-connected, and hence has minimum degree at least $3$.

\begin{thm}{\rm (Ho, Lin, Tan, Hsu and Hsu \cite{4})}
Let $\phi(n,\delta)$ be the maximum sizes of a nonhamiltonian-connected graph of order $n$ and minimum degree at least $\delta$. Then
$$
\phi(n,\delta)=\begin{cases}
\binom{n-\delta+1}{2}+\delta(\delta-1)\quad {\rm if}\,\,n\,\,\,{\rm is}\,\,{\rm odd}\,\,{\rm and}\,\,3\le \delta \le \frac{n+9}{6},\\
\qquad\qquad\qquad\qquad\,\,\,\,\, {\rm or}\,\, n\,\,{\rm is}\,\,{\rm even}\,\,{\rm and}\,\,3\le \delta \le \frac{n+6}{6},\\
\binom{n-\lfloor\frac{n}{2}\rfloor+1}{2}+\lfloor \frac{n}{2}\rfloor(\lfloor \frac{n}{2}\rfloor-1)\quad{\rm if}\,\,n\,\,\,{\rm is}\,\,{\rm odd}\,\,{\rm and}\,\,\frac{n+9}{6}< \delta \le \frac{n-1}{2},\\
\qquad\qquad\qquad\qquad\qquad\quad\,\,\,\,{\rm or}\,\, n\,\,{\rm is}\,\,{\rm even}\,\,{\rm and}\,\,\frac{n+6}{6}< \delta \le \frac{n}{2}.
\end{cases}
$$
\end{thm}

Note that if $n$ is odd and $\frac{n+9}{6}< \delta < \frac{n-1}{2}$ or $n$ is even and  $\frac{n+6}{6}< \delta < \frac{n}{2}$, the maximum size of Theorem 6 is not attained by a graph of minimum degree $\delta.$ For example, if $n=16$ and $\delta=4$, the maximum size $\phi(16,4)=92$ is attained by a unique graph of minimum degree $8,$ not $4.$ In this paper we obtain more precise information by determining the maximum  number of $s$-cliques in  a nonhamiltonian-connected  graph of order $n$ and minimum degree $\delta.$ As a corollary, we determine the maximum size of a nonhamiltonian-connected graph with given order and minimum degree, as well as the extremal graphs, from which Theorem 5 and Theorem 6 can be deduced. Before presenting the main theorem, we need the following notations.

\begin{nota}
Fix $3\le \delta\le \lfloor \frac{n}{2}\rfloor.$ Let $F(n,\delta)=K_{\delta}\vee(K_{n-2\delta+1}+\overline{K_{\delta-1}})$. Denote by $f_s(n,\delta)$ the number of $s$-cliques in $F(n,\delta)$; more precisely,
$$f_s(n,\delta)=\binom{n-\delta+1}{s}+(\delta-1)\binom{\delta}{s-1}.$$
 \end{nota}

\begin{nota}
For $t=\lfloor n/2\rfloor$, let $G(n,\delta)$ denote the graph obtained from $K_t\vee(K_{n-2t+1}+\overline{K_{t-1}})$  by deleting $t-\delta$ edges that are incident to one common vertex in $\overline{K_{t-1}}$. Denote by $g_s(n,\delta)$ the number of $s$-cliques in $G(n,\delta)$; more precisely,
\begin{align*}
  g_s(n,\delta)=\binom{n-t+1}{s}+(t-2)\binom{t}{s-1}+\binom{\delta}{s-1}.
\end{align*}
\end{nota}

It can be checked that $f_s(n,\delta)$ and $g_s(n,\delta)$ are convex functions.

Note that if the minimum degree $\delta=1,$ the maximum size of a  nonhamiltonian-connected graph of order $n$ is $\binom{n-1}{2}+1$, and if $\delta=2,$ the maximum size is  $\binom{n-1}{2}+2$. Therefore, in the following, we will assume that $\delta\ge 3$. The following are the main results of this paper.

\begin{thm}\label{th1}
 Let $\varphi_s(n,\delta)$ denote the maximum number of $s$-cliques in a nonhamiltonian-connected graph of order $n$ and minimum degree $\delta$ where $3\le \delta\le\lfloor \frac{n}{2}\rfloor $.
 Then
$$
\varphi_s(n,\delta)= {\rm max}\{f_s(n,\delta),g_s(n,\delta)\}.
$$
\end{thm}

This theorem is sharp with extremal graphs $F(n,\delta)$ and $G(n,\delta)$.

\begin{thm}\label{th2}
 Let $\varphi(n,\delta)$ denote the maximum size of a nonhamiltonian-connected graph of order $n$ and minimum degree $\delta$ where $3\le \delta\le\lfloor \frac{n}{2}\rfloor $. Then
$$
\varphi(n,\delta)=\begin{cases}
\binom{n-\delta+1}{2}+\delta(\delta-1)\,\,\ {\rm if}\,\,n\,\,{\rm is}\,\,{\rm odd}\,\,{\rm and}\,\,3\le \delta \le \frac{n+13}{6}\,\,{\rm or}\,\,
n\,\,{\rm is}\,\,{\rm even}\,\,{\rm and}\,\,3\le \delta \le \frac{n+10}{6},\\
\frac{3n^2-8n+13}{8}+\delta\quad{\rm if}\,\,n\,\,{\rm is}\,\,{\rm odd}\,\,{\rm and}\,\,\frac{n+13}{6}\le \delta\le\frac{n-1}{2},\\
\frac{3n^2-6n}{8}+\delta\quad {\rm if}\,\,n\,\,{\rm is}\,\,{\rm even}\,\,{\rm and}\,\,\frac{n+10}{6}\le \delta\le\frac{n}{2}.
\end{cases}
$$
If $\delta=\frac{n+13}{6}$  or $\delta=\frac{n+10}{6}$, then $\varphi(n,\delta)$ is attained by a graph $Q$ if and only if $Q=F(n,\delta)$ or $Q=G(n,\delta).$
If $n$ is odd and $3\le\delta\le\frac{n+11}{6}$  or $n$ is even and $3\le\delta\le\frac{n+8}{6}$, then $\varphi(n,\delta)$ is attained by a graph $Q$ if and only if $Q=F(n,\delta).$
If $n$ is odd and $\frac{n+15}{6}\le \delta\le\frac{n-1}{2}$ or $n$ is even and $\frac{n+12}{6}\le \delta\le\frac{n}{2},$ then $\varphi(n,\delta)$ is attained by a graph $Q$ if and only if $Q=G(n,\delta).$
\end{thm}

The rest of this paper is organized as follows. In the next section we collect various known results that will be used in later arguments and then give a proof of Theorem \ref{th1}. In Section 3 we present a proof of Theorem \ref{th2} and some corollaries.

\section{Proof of Theorem \ref{th1}}

We will need the following lemmas.

\begin{lem}\label{lem1}{\rm (Ore \cite{2})}
Let $n\ge3$ and $G$ an $n$-vertex graph. If $d(u)+d(v)\ge n+1$ for every pair $u, v\in V(G)$ with $u\nleftrightarrow v$, then $G$ is hamiltonian-connected. In particular, if $\delta(G)\ge(n + 1)/2$, then $G$ is hamiltonian-connected.
\end{lem}

\begin{lem}\label{lem2}{\rm (Bondy and Chv\'{a}tal \cite{7})}
Suppose that there are two vertices $u$ and $v$ in $V(G)$ satisfying $d(u)+d(v)\ge n+1$ and $u\nleftrightarrow v.$ Then $G+uv$ is hamiltonian-connected if and only if $G$ is hamiltonian-connected.
\end{lem}

To prove our results, we also need a definition from Kopylov.

\begin{defi}{\rm ($t$-disintegration of a graph)}
 Let $G$ be a graph and $t$ be a positive integer. Delete all vertices of degree at most $t$ from $G$; for the resulting graph $G'$, we again delete all vertices of degree at most $t$ from $G'$. Iterating this process until we finally obtain a graph, denoted by $D(G; t)$, such that either $D(G; t)$ is a null graph or $\delta(D(G; t))\geq t+1.$ The graph $D(G; t)$ is called the {\it $(t+1)$-core of $G.$}
\end{defi}

For $s\geq 2$, let $N_s(G)$ denote the number of $s$-cliques in $G$. Now we are ready to prove Theorem \ref{th1}.

\noindent{\bf Proof of Theorem \ref{th1}.} Let $P$ be a hamiltonian path in graph $W$ with endpoints $u$ and $v$ and let $S$ be a vertex set in $W$ such that $u,v\in S.$ Then the number of components in the graph $W-S$ is at most $|S|-1.$  Using this fact, we can verify that the graphs $F(n,\delta)$ and $G(n,\delta)$ defined in Notation 1 and Notation 2, respectively, are nonhamiltonian-connected graphs of order $n$ and minimum degree $\delta.$ The number of copies of $s$-cliques in $F(n,\delta)$ ($G(n,\delta)$ respectively) is $f_s(n,\delta)$ ($g_s(n,\delta)$ respectively).

Let $G$ be an $n$-vertex nonhamiltonian-connected graph with $\delta(G)=\delta.$ Let $w$ be a vertex of $G$ with minimum degree $\delta.$ If there exist two vertices $u,v\in V(G)\backslash\{w\}$ such that $u\nleftrightarrow v$ and $d_G(u)+d_G(v)\ge n+1,$ we denote by $G_1$ the graph $G+uv.$ For the graph $G_1,$ we again choose $u_1,v_1\in V(G_1)\backslash\{w\}$ with $u_1\nleftrightarrow v_1, d_{G_1}(u)+d_{G_1}(v)\ge n+1,$ and denote by $G_2$ the graph $G_1+u_1v_1.$ Iterating this process until we finally obtain a graph, denoted by $Q,$ such that for any $x,y\in V(Q)\backslash \{w\}$ and $x\nleftrightarrow y,$ we have $d_Q(x)+d_Q(y)\le n.$ We will repeatedly use this condition without mentioning it possibly. Obviously, $\delta(Q)=\delta$ and the graph $Q$ is a nonhamiltonian-connected graph by Lemma \ref{lem2}.

Let $t=\lfloor\frac{n}{2}\rfloor.$ Denote by $D=D(Q;t)$ the $(t+1)$-core of $Q$, i.e. the resulting graph after applying $t$-disintegration to $Q.$ We distinguish two cases.

{\bf Case 1.} $D$ is a null graph. In the $t$-disintegration process, put $Q_0=Q$ and $Q_{i+1}=Q_i-x_i,0\leq i\leq n-1$ where $x_i$ is a vertex of degree at most $t$ in $Q_i$.
Since $\delta(Q)\leq t$ (otherwise $\delta(Q)\geq t+1$ and no vertex could be deleted), without loss of generality, we can take $x_0$ to be $w$, the vertex with minimum degree.
Note that once the vertex $x_i$ is deleted, we have deleted at most $\binom{d_{Q_i}(x)}{s-1}$ copies of $K_s$. By the definition of $t$-disintegration, we have $d_{Q_i}(x_i)\leq t,1\leq i\leq n-t-1.$  For the last $t$ vertices, the number of $K_s$ is at most $\binom{t}{s}$. Thus
$$
  N_s(Q)\leq\binom{\delta}{s-1}+(n-t-1)\binom{t}{s-1}+\binom{t}{s}\le g_s(n,\delta).
$$

{\bf Case 2.} $D$ is not a null graph. Let $d=|D|.$ We claim that $V(D)$ is a clique and $\delta\le n+1-d$.

For all $u, v\in V(D)$, we have $d_D(u),\ d_D(v)\geq t+1$. Since every nonadjacent pair of vertices has degree sum at most $n$ in $Q$ and $d_Q(u)+d_Q(v)\geq d_D(u)+d_D(v)\ge 2t+2\ge n+1$, we have $u$ and $v$ are adjacent in $Q$, i.e., $V(D)$ is a clique.

We next prove $\delta\le n+1-d.$ Suppose on the contrary that $d\geq n+2-\delta$. Then $d_D(u)\geq d-1\geq n+1-\delta$ for all $u\in V(D)$. Since $V(D)$ is a clique and $d_D(u)\geq t+1$ for all $u\in V(D),$ we have $d\geq t+2.$ Thus, every vertex in $ V(Q)\backslash V(D)$ is not adjacent to at least two vertices in $D.$ Let $x\in V(Q)\backslash V(D)$ and $y\in V(D)$ not adjacent to $x$. Note that $w\in V(Q)\setminus V(D).$ We distinguish two cases. If $V(Q)\setminus V(D)=\{w\},$ we have $x=w$ and $|D|=n-1.$ Then $Q$ is hamiltonian-connected as $D$ is a complete graph, a contradiction. If $V(Q)\setminus V(D)\ne\{w\},$ we may assume $x\ne w.$ Note that $d_Q(x)\geq \delta,$ we have $d_Q(x)+d_Q(y)\geq\delta+n+1-\delta=n+1$, a contradiction due to the structure of graph $Q$. Thus $d\le n+1-\delta,$ i.e., $\delta\le n+1-d$.

Let $D'$ be the $(n+2-d)$-core of $Q,$ i.e. the resulting graph after applying $(n+1-d)$-disintegration to $Q.$ Since $d\geq t+2$, we obtain $n+1-d\leq n+1-t-2\le t.$ Therefore, $D\subseteq D'$. There are two cases.

(a) If $D'=D$, then $|D'|=|D|=d.$ By the definition of $(n+1-d)$-disintegration, we have
\begin{align*}
N_s(Q)
&\leq \binom{\delta}{s-1}+(n-d-1)\binom{n+1-d}{s-1}+\binom{d}{s}\\
&=\binom{\delta}{s-1}+\lambda_s(n,n+1-d)\\
&\leq {\rm max}\{ f_s(n,\delta),g_s(n,\delta)\},
\end{align*}
where $\lambda_s(n,x)=(x-2)\binom{x}{s-1}+\binom{n+1-x}{s}.$ The third inequality follows from the condition $\delta\leq n+1-d \leq t$ and that the function $\lambda_s(n,x)$ is convex for $x\in[\delta,t].$

(b) Otherwise $D'\neq D.$ Let $u\in V(D')\backslash V(D).$ Since $d\ge t+2,$ we deduce that $u$ is not adjacent to at least two vertices in $D$. We choose one of the vertices and denote it by $v$. Then $d_Q(u)+d_Q(v)\ge n+2-d+d-1\ge n+1.$ Since every nonadjacent pair of vertices has degree sum at most $n,$ we obtain a contradiction. Hence the proof of Theorem \ref{th1} is complete.\hfill $\Box$

\section{Proof of Theorem \ref{th2} and some corollaries}

In order to present extremal graphs in Theorem \ref{th2}, we need the following lemma.

\begin{lem}\label{lem3}{\rm (Lick \cite{5})}
Let $G$ be a graph with degree sequence $d_1\le d_2\le\cdots\le d_n$ where $n\ge 3.$ If there is no integer $i$ with $2\le i\le n/2$ such that $d_{i-1}\le i$ and $d_{n-i}\le n-i,$ then $G$ is hamiltonian-connected.
\end{lem}

\noindent{\bf Proof of Theorem \ref{th2}.}
The expression of function $\varphi(n,\delta)$ can be easily obtained from Theorem \ref{th1}. Now we determine the extremal graphs. Since the proof when $n$ is even is similar to the proof when $n$ is odd, we give only the proof of the latter.

Let $Q$ be a nonhamiltonian-connected graph of order $n$ and minimum degree $\delta.$ The degree sequence of $Q$ is $d_1\le d_2\le\cdots\le d_n$ where $n\ge 3.$  By Lemma \ref{lem3}, there exists an $i$ with $i\le(n-1)/2$ such that $d_{i-1}\le i$ and $d_{n-i}\le n-i.$ We shall give extremal graphs with $d_{i-1}=i$ and $d_{n-i}=n-i$. Now assume that each vertex degree attains its possible maximum. More precisely, let us consider the following degree sequence:
$$
\delta,\,\,\, \underbrace{i,\ldots,i}_{i-2},\,\,\, \underbrace{n-i,\ldots,n-i}_{n-2i+1},\,\,\, \underbrace{n-2,\ldots,n-2}_{i-\delta},\,\,\, \underbrace{n-1,\ldots,n-1}_{\delta}.
$$
This degree sequence is graphical and the sum of all vertex degrees is $n^2-(2i-1)n+3i^2-5i+2\delta.$

In the case $3\le\delta\le (n+13)/6,$ we have $n\ge 6\delta-13.$ Suppose that $Q$ has size $\binom{n-\delta+1}{2}+\delta(\delta-1).$  We have
$$
e(Q)=\binom{n-\delta+1}{2}+\delta(\delta-1)\le \frac{n^2-(2i-1)n+3i^2-5i+2\delta}{2}, \eqno (1)
$$
where the inequality is equivalent to $(i-\delta)(2n-3i-3\delta+5)\le 0.$ Since $i\ge d_i\ge \delta(Q)=\delta,$ we obtain $i=\delta$ or $i>\delta$ and $n\le (3i+3\delta-1)/2.$

If $i=\delta,$ equality holds in (1) and hence the degree sequence of $Q$ is
$$
\underbrace{\delta,\ldots,\delta}_{\delta-1},\, \underbrace{n-\delta,\ldots,n-\delta}_{n-2\delta+1},\, \underbrace{n-1,\ldots,n-1}_{\delta},
$$
implying that $Q=K_\delta\vee(K_{n-2\delta+1}+\overline{K_{\delta-1}}).$

Now suppose $i> \delta.$ Then we have $n\le (3i+3\delta-5)/2.$ If $i\le (n-3)/2,$ then $n\le 6\delta-19,$ contradicting to our assumption $n\ge 6\delta-13.$ Thus $i=(n-1)/2.$
We have $n\leq 6\delta-13.$ Note that $n\ge 6\delta-13,$ implying that $n=6\delta-13.$  Hence the degree sequence of $Q$ is
$$
\delta,\, \underbrace{\frac{n-1}{2},\ldots,\frac{n-1}{2}}_{\frac{n-1}{2}-2},\, \underbrace{\frac{n+1}{2},\frac{n+1}{2}}_{2}, \underbrace{n-2,\ldots,n-2}_{\frac{n-1}{2}-\delta},\, \underbrace{n-1,\ldots,n-1}_{\delta},
$$
implying that $Q=G(n,\delta).$

Now we prove Theorem \ref{th2} for $\frac{n+15}{6}\le \delta\le\frac{n-1}{2}.$ Obviously, $2\delta+1\le n\le 6\delta-15.$ Suppose that $Q$ has size $\frac{3n^2-8n+13}{8}+\delta.$ We have
$$
e(Q)=\frac{3n^2-8n+13}{8}+\delta\le \frac{n^2-(2i-1)n+3i^2-5i+2\delta}{2}, \eqno (2)
$$
where the inequality is equivalent to $0\le (6i-13-n)(2i+1-n).$ Since $i\ge d_i\ge \delta(Q)=\delta$ and $n\le 6\delta-15.$ Then $n\leq 6i-15.$  The inequality (2) is equivalent to $i\ge (n-1)/2.$ Note that $(n-1)/2\ge i$. We have $i=(n-1)/2.$ Hence the degree sequence of $Q$ is
$$
\delta,\, \underbrace{\frac{n-1}{2},\ldots,\frac{n-1}{2}}_{\frac{n-1}{2}-2},\, \underbrace{\frac{n+1}{2},\frac{n+1}{2}}_{2}, \underbrace{n-2,\ldots,n-2}_{\frac{n-1}{2}-\delta},\, \underbrace{n-1,\ldots,n-1}_{\delta},
$$
implying that $Q=G(n,\delta).$ This completes the proof. \hfill $\Box$

The following corollary follows from Theorem \ref{th2} immediately.

\begin{cor} Let $\varphi(n,\delta)$ be defined as in Theorem \ref{th2}. If $G$ is a graph of order $n$ and minimum degree $\delta$ with size greater than $\varphi(n,\delta),$ then $G$ is hamiltonian-connected.
\end{cor}

In 1995 Ota \cite{17} determined the maximum size of a $k$-connected nonhamiltonian graph of order $n.$ Therefore, it is natural to consider the corresponding problem for a nonhamiltonian-connected graph.

\begin{cor}
Let $G$ be a $k$-connected ($k\geq 3$) nonhamiltonian-connected graph of order $n$.  Then
\begin{equation*}
  N_s(G)\leq
  \mathop{{\rm max}}\limits_{k\leq \delta\leq \lfloor n/2\rfloor}
  \left\{{\rm max}\{f_s(n,\delta),g_s(n,\delta)\}\right\}
  ={\rm max}\{f_s(n,k),\,g_s(n,\lfloor n/2\rfloor) \}.
\end{equation*}
\end{cor}

\noindent{\bf Acknowledgment.} This research was supported by the NSFC grant 12271170 and Science and Technology Commission of Shanghai Municipality (STCSM) grant 22DZ2229014.


\begin{thebibliography}{99}
\bibitem{10}J.A. Bondy, Variations on the hamiltonian theme, Canad. Math. Bull. 15 (1972)  57-62.
\bibitem{7}J.A. Bondy, V. Chv\'{a}tal, A method in graph theory, Discrete Math. 15 (1976) 111-135.
\bibitem{6}P. Erd\H{o}s, Remarks on a paper of P\'{o}sa, Magyar Tud. Akad. Math. Kutat\'{o}. Int. K\"{o}zl. 7 (1962) 227-229.
\bibitem{3}Z. F\"{u}redi, A. Kostochka, R. Luo, Extensions of a theorem of Erd\H{o}s on nonhamiltonian graphs, J. Graph Theory 89 (2018) 176-193.
\bibitem{11}M.R. Garey, D.S. Johnson, L. Stockmeyer, Some simplified NP-complete graph problems, Theoret. Comput. Sci. 1 (1976) 237-267.
\bibitem{4}T.Y. Ho, C.K. Lin, J.J.M. Tan, D.F. Hsu, L.H. Hsu, On the extremal number of edges in hamiltonian connected graphs, Appl. Math. Lett. 23 (2010) 26-29.
\bibitem{12}T.Y. Ho, C.K. Lin, J.J.M. Tan, D.F. Hsu, L.H. Hsu, On the Extremal Number of Edges in Hamiltonian Graphs, J. Inf. Sci. Eng. 27 (2011) 1659-1665.
\bibitem{19}Z. Khan, L-T Yuan, A note on the 2-power of Hamilton cycles, Discrete Math. 345 (2022) Paper No.112908.
\bibitem{14}B.L. Li, B. Ning, Spectral analogues of Erd\H{o}s and Moon-Moser's theorems on Hamiltoncycles, Linear Multilinear Algebra 64 (2016) 2252-2269.
\bibitem{15}B.L. Li, B. Ning, X. Peng, Extremal problems on the Hamiltonicity of claw-free graphs, Discrete Math. 341 (2018) 2774-2788.
\bibitem{5}D.R. Lick, A sufficient condition for hamiltonian connectedness, J. Combin. Theory 8 (1970) 444-445.
\bibitem{20}X. Liu, Z. Ryj\'{a}\v{c}ek, P. Vr\'{a}na, L.M. Xiong and X.J. Yang, Hamilton-connected \{claw,net\}-free graphs, II, J. Graph Theory 103 (2023) 119-138.
\bibitem{17} K. Ota, Cycles through prescribed vertices with large degree sum, Discrete Math. 145 (1995) 201-210.
\bibitem{2}O. Ore, Hamiltonian connected graphs, J. Math. Pures Appl. 42 (1963) 21-27.
\bibitem{8}O. Ore, Arc coverings of graphs, Ann. Mat. Pura Appl. 55 (1961) 315-321.
\bibitem{1}E.F. Schmeichel, S.L. Hakimi, Pancyclic Graphs and a Conjecture of Bondy and Chv\'{a}tal, J. Combin. Theory Ser. B 17 (1974) 22-34.
\bibitem{9}D.B. West, Introduction to Graph Theory, Prentice Hall, Inc., 1996.
\bibitem{18}X. Zhan, The maximum degree of a minimally hamiltonian-connected graph, Discrete Math. 345 (2022) Paper No. 113159.
\bibitem{16}L. Zhang, The maximum number of cliques in graphs with prescribed order, circumference and minimum degree, European J. Combin. 112 (2023) Paper No. 103728.
\bibitem{13}Q.N. Zhou, H. Broersma, L.G. Wang, Y. Lu, On sufficient spectral radius conditions for hamiltonicity of $k$-connected graphs, Linear Algebra Appl. 604 (2020) 129-145.
\end{thebibliography}
\end{document}